\documentclass[english]{amsart}

\usepackage[T1]{fontenc}
\usepackage[utf8]{inputenc}
\usepackage{microtype, amsmath, amsfonts, amssymb, amsthm, amscd, braket, bm, bookmark, graphicx, mathrsfs, mathtools, mathdots, nicematrix,tikz, caption, subfig, import, centernot,enumitem,url,ragged2e,tikz-cd, hyperref}
\usepackage[style=alphabetic, maxnames=5, backend=biber,sorting=nyt]{biblatex}
	\renewbibmacro{in:}{\ifentrytype{article}{}{\printtext{\bibstring{in}\intitlepunct}}}
    
\usetikzlibrary{automata,positioning,calc,intersections,through,backgrounds,patterns,fit,external}

\theoremstyle{definition}
\newtheorem{defn}{Definition}[section]
\newtheorem*{defn*}{Definition}

\newtheorem{quest}[defn]{Question}

\newtheorem*{fact*}{Fact}

\theoremstyle{plain}
\newtheorem{thm}[defn]{Theorem}
\newtheorem*{thm*}{Theorem}

\newtheorem{cor}[defn]{Corollary}
\newtheorem{lem}[defn]{Lemma}
\newtheorem*{conj*}{Conjecture}

\theoremstyle{remark}
\newtheorem{rmrk}[defn]{Remark}
\newtheorem*{rmrk*}{Remark}

\newenvironment{Proof}
{\par{\noindent\bf Proof:}}
{\hfill$\blacksquare$}

\newcommand{\norm}[1]{\left\lVert #1 \right\rVert}
\newcommand{\SizedSet}[1]{\left\{\, #1 \,\right\}}

\providecommand{\given}{}

\newcommand{\SetSymbol}[1][]{
	\nonscript\:\mathord{:}
	\allowbreak
	\nonscript\:
	\mathopen{}}

	\renewcommand{\given}{\SetSymbol[\delimsize]}

\DeclarePairedDelimiter{\abs}{\lvert}{\rvert}
\DeclareMathOperator{\cotan}{cotan}  
\DeclareMathOperator{\Leb}{Leb}  

\bibliography{biblio}

\begin{document}

\title[Triangular Billiards of Weakly Exponential Growth]{Explicit Families and distribution of Triangular Billiards of Weakly Exponential Growth}

\author[I. Mamsurova]{Irina Mamsurova}
\address{Aix Marseille Université, CNRS, Centrale Marseille, I2M, UMR 7373, 13453 Marseille, France.}
\email{\href{mailto:irina.mamsurova@univ-amu.fr}{irina.mamsurova@univ-amu.fr}}

\begin{abstract}
Recently, it has been shown that the combinatorial complexity function $N_c(n)$ of a typical triangular billiard has weakly exponential growth, i.e., for almost any triangle and any $\varepsilon>0$ there is a constant $C$ such that $N_c(n)<Ce^{n^\varepsilon}$. We give the first example, an infinite family of explicitly described triangles with weakly exponential complexity growth. The simplest example in this family is a right triangle with integer side lengths. Moreover, we find one-parameter families of triangular tables with weakly exponential growth for generic parameters. 
\end{abstract}

\maketitle

\section{Introduction}

\subsection{Complexity of polygonal billiards.} We study the billiard motion in a polygon $Q$. A point mass moves in $Q$ along a straight line until it hits a side, i.e. a point of $\partial Q$. Then it bounces from the side according to the law of geometric optics ``the angle of incidence is equal to the angle of reflection'' and continues the motion along the new line. If the point hits a vertex, then its future motion is not defined.

The question of the complexity of this dynamical system can be formulated in terms of generalized diagonals. The \textit{generalized diagonal} is an oriented billiard trajectory that starts and terminates at a vertex (not necessarily the same). Generalized diagonals can be viewed as ``highly degenerate'' trajectories. By the combinatorics of an orbit we mean the sequence of sides that it hits. The combinatorics of two orbits of a bounded length, starting at a vertex, coincide if and only if between them there is no generalized diagonal of a smaller length starting at the same vertex. We can reformulate it in terms of the phase space. The parameters of the phase space are the point of $\partial Q$ and the angle between the trajectory and $\partial Q$. The phase space is divided into $n$-cells by ``degenerate'' trajectories  having at most $n$ reflections. Combinatorics of the first $n$ reflections of two orbits corresponding to points lying in the same $n$-cell of a phase space are the same (see, for example, \cite{Kat1} or \cite{CHT}). The cells are separated by segments corresponding to orbits, starting or terminating at a vertex; the points of intersection of these segments correspond to generalized diagonals.

The combinatorial length $\ell_c$ of a generalized diagonal is defined by the number of reflections of sides that the generalized diagonal has (note that in some texts the combinatorial length is defined as the number of links in a generalized diagonal, as in \cite{CHT}). Its geometric length $\ell_g$ is a sum of euclidian lengths of its links. 

\begin{defn}
   The \textit{combinatorial complexity function} $N_c(n)$  corresponding to a billiard in a polygon is the number of generalized diagonals of combinatorial length at most $n\in\mathbb{N}$: 
 \[\ell_c(\gamma)=\#\text{ reflections of sides of a trajectory }\gamma,\]
 \[N_c(n)=\#\SizedSet{\gamma \given \ell_c(\gamma)\leq n}.\]
 
The \textit{geometric complexity function} $N_c(t)$  corresponding to a billiard in a polygon is the number of generalized diagonals of combinatorial length at most $t\in\mathbb{R}_+$: 
\[\ell_g(\gamma)=\sum_{v \text{ is a link of }\gamma}\text{ euclidian length of }v,\] 
 \[N_g(t)=\#\SizedSet{\gamma \given \ell_g(\gamma)\leq t}.\] 
\end{defn} 
The results for combinatorial and geometric complexity are related by the existence of $C>1$ such that $C^{-1}\leq\frac{N_c(t)}{N_g(t)}\leq C$.

An alternative way to describe the complexity of a billiard is the billiard language. Label the sides of $Q$ by distinct letters of an alphabet $\mathcal{A}$ whose cardinality is equal to the number of sides of $Q$. An orbit of combinatorial length $n$ codes a word of length $n$ in the alphabet $\mathcal{A}$, defined by the sequence of sides it hits. We call the resulting language \textit{the billiard language}; it consists of words describing all possible combinatorics of orbits. Two billiard trajectories code words with the same prefix of length $n$ if and only if the corresponding points in the phase space belong to the same $n$-cell. The complexity $p(n)$ of the associated billiard language is the number of words of length $n$.

The following result relates two quantities that measure the complexity of a billiard. It was proved for convex polygons by J. Cassaigne, P. Hubert, and S. Troubetzkoy in \cite{CHT}. T. Krueger, A. Nogueira, and S. Troubetzkoy observed in \cite{KNT} that the same proof also works for non-convex simply connected polygons.

\begin{thm}[J. Cassaigne, P. Hubert, S. Troubetzkoy, \cite{CHT}; T. Krueger, A. Nogueira, S. Troubetzkoy, \cite{KNT}]\label{CHT}
    For any simply connected polygon $Q$
    \[p(n)=\sum_{j=0}^{n-2}N_c(j), \]
where the vertices of the polygon are counted as generalized diagonals of combinatorial length $0$
and the sides of $Q$ are not counted as generalized diagonals.
\end{thm}
A polygon is called rational if its inner angles are in $\pi\mathbb{Q}$. In \cite{Ma1} and \cite{Ma2} H. Masur showed that for any rational polygon the complexity function $N_g(n)$ has quadratic upper and lower bounds: there are constants $C_1, C_2>0$ such that $C_1\cdot t^2<N_g(t)<C_2\cdot t^2$. It follows that there are quadratic bounds for combinatorial complexity $N_c(n)$ and (by Theorem \ref{CHT}) cubic bounds for language complexity $p(n)$. W. A. Veech \cite{V1}, \cite{V2} introduced a class of rational polygons such that the quantity $N_g(t)/t^2$ admits a limit as $t$ tends to infinity. J. Athreya, P. Hubert, and S. Troubetzkoy show that for polygons of this class $N_c(n)/n^2$ and $p(n)/n$ admits a limit as $n$ tends to infinity and give the value of these limits for regular $k$-gons \cite{AHT}.   

Much less is known about irrational polygons. The difficulty in this case is that we cannot use the machinery of the associated compact translation surface. A. Katok \cite{Kat1} proved the subexponential estimate for the complexity function: for any polygon $\lim_{t\to+\infty}\frac{\ln(N_g(t))}{t}=0$.

The gap between the known non-explicit sub-exponential
estimate and the quadratic growth known for the rational case motivated the third question in A. Katok's list ``Five most resistant problems in dynamics''.

\begin{quest}[A. Katok, \cite{Kat2}]
    Orbit growth, complexity growth, existence of periodic
orbits and ergodicity of billiards in general (not rational) polygons.
\end{quest}

T. Krueger, A. Nogueira and S. Troubetzkoy in \cite{KNT} showed that for billiard in a typical polygon the geometric and combinatorial complexity $N_g(t)$ and $N_c(n)$ grow quadratically and the language complexity $p(n)$ grows cubically along certain subsequences. Recently, D. Scheglov provided the polynomial lower and weakly exponential upper bounds of a combinatorial complexity function for typical triangles in \cite{Sch1}, \cite{Sch2}. Both bounds reflect the typical behavior of complexity functions: the lower bound is given for a generic right triangle (\cite{Sch1}), and the upper bound is given for a generic triangle (\cite{Sch2}).  We focus on the weakly exponential growth result. Weakly exponential growth means precisely that for every $\varepsilon>0$ there is a constant $C>0$ such that $N_c(n)<Ce^{n^\varepsilon}$. 

\begin{thm}[Weakly exponential growth, D. Scheglov, \cite{Sch2}]\label{weaklyexp} For a typical triangle and any $\varepsilon>0$ there is a constant $C>0$ such that $N_c(n)<Ce^{n^{\varepsilon}}$.
\end{thm}

Using Theorem \ref{CHT}, we obtain an immediate corollary for the complexity of a billiard language.  
\begin{cor} For a typical triangle and any $\varepsilon>0$ there is a constant $C>0$ such that $p(n)<Ce^{n^{\varepsilon}}$.
\end{cor}

The proof of Theorem \ref{weaklyexp} is based on a weak Diophantine condition for almost all elements of $SO(3)$ (or, equivalently, $SU(2)$). The weak Diophantine condition for Haar-almost all elements of $SO(3)$ was proven by V. Kaloshin and I. Rodnianski in \cite{KR}. Precisely, for almost all $A,B\in SO(3)$ there is a constant $D>0$ such that for every $n\in\mathbb{N}$ we have $\norm{W_n(A,B)-E}\ge D^{-n^2}$, where $E$ denotes the identity in $SO(3)$ and $W_n(A,B)$ denotes a reduced word in $A^{\pm 1},B^{\pm 1}$ of length $n$. The strong Diophantine condition of elements of $SU(2)$ with algebraic angles was proven by A. Gamburd, D. Jacobson and P. Sarnak (\cite{GJS}, Proposition 4.3): for any $k\in\mathbb{N}$ and any $g_1,\dots,g_k\in SU(2)$ with algebraic entries there is a constant $D>0$ such that for every $n\in\mathbb{N}$ we have $\norm{W_n(g_1,\dots,g_k)-E}\ge D^{-n}$. Again, $E$ denotes the identity and $W_n(g_1,\dots,g_k)$ denotes a reduced word in $g_1^{\pm 1},\dots,g_k^{\pm 1}$ of length $n$. For elements of $SO(3)$ (equivalently, $SU(2)$) with general entries the question of strong Diophantine condition posed in \cite{GJS} remains open. 

\subsection{Main results.} Even though the previous results hold for a typical triangle, no examples of irrational billiards with weakly exponential growth were known until now. However, the work \cite{Sch2} provides a criterion for weakly exponential growth for a triangle. This criterion allows one to construct  an infinite family of examples.

The main theorem of this paper shows that an explicitly described infinite family of triangles gives rise to billiards with weakly exponential growth. 

\begin{defn}
    We call an angle $\theta$ \textit{trigonometrically algebraic} (respectively, \textit{trigonometrically rational}) if $\sin\theta$, $\cos\theta$ are algebraic (respectively, rational) real  numbers.
\end{defn}

\begin{thm}\label{main_thm}
    Let the angles $\alpha,\beta$ be trigonometrically algebraic. Then the growth of the combinatorial complexity function in the $(\alpha,\beta,\pi-\alpha-\beta)$ triangle is weakly exponential.
\end{thm}

Again, Theorem \ref{CHT} yields the following corollary.
\begin{cor}
    Let the angles $\alpha,\beta$ be trigonometrically algebraic. Let $p(n)$ denote the complexity of the billiard language in the $(\alpha,\beta,\pi-\alpha-\beta)$ triangle. Then for every $\varepsilon>0$ there is a constant $C$ such that $p(n)<Ce^{n^\varepsilon}$.
\end{cor}

The simplest example of a triangle with known weakly exponential growth is a right triangle with integer side lengths. It corresponds to the case of trigonometrically rational angles, for which the proof of weakly exponential growth is very simple. The proof for trigonometrically algebraic angles is based on the Liouville inequality. Theorem \ref{main_thm} is, in some sense, analogous to the strong Diophantine condition for elements of $SU(2)$ proved in \cite{GJS}. The method presented by D. Scheglov in \cite{Sch2} and developed in the following text permits one to pass from the geometry of a billiard table to the arithmetic properties of the triginometric values of its angles.

To the best of our knowledge, the measure distribution of triangles of weakly exponential growth was until now unknown. For example, is it true that almost all triangles with the smallest angle equal to $\pi/4$ the growth is weakly exponential? We  obtain a result about the distribution of weakly-exponential billiard tables in the set of triangles with one fixed angle.

\begin{thm}\label{repart}
    Let $\alpha$ be a fixed trigonometrically algebraic angle.  Then for almost all $\beta$ the complexity function of a billiard in the triangle $(\alpha,\beta,\pi-\alpha-\beta)$ is weakly exponential.
\end{thm} 

\subsection{Structure of the paper.} We recall some results of \cite{Sch2} and formulate explicitly the criterion for weakly exponential growth in Section \ref{sec2}. Then in Corollary \ref{rational} we give a simple example of an explicitly described (infinite) family of triangles that satisfies this criterion. Generalizing the arguments, we prove Theorem \ref{main_thm}. Section \ref{sec3} is dedicated to the distribution (in measure sense) of weakly exponential triangular tables with a fixed trigonometrically algebraic angle; again, we first describe the simple example of a right triangular table and then move on to the general case, proving Theorem \ref{repart}.

\section*{Acknowledgements} 
This work started during my Master thesis under the supervision of Pascal Hubert in Institut de Mathématiques de Marseille; Pascal's guidance and enthusiasm made the appearance of this  text possible. I am grateful to Serge Troubetzkoy, Tyll Krueger, Giovanni Forni and Carlos Matheus as well as to the participants of Rauzy seminar and the eighth winter school in geometry and dynamics in Aussois for being the first audiences of the presented results. I also sincerely thank Yann Bugeaud for suggesting the use of the Liouville inequality and Thibaut Misme for explaining its general form. 

I was partially supported by Fondation Sciences mathématiques de Paris and Fondation de l’École normale supérieure. 
\section{The criterion for weakly exponential growth}\label{sec2}

\begin{figure}
    \centering
    \begin{tikzpicture}[scale=1.0, line join=round, line cap=round]
 
\coordinate (A) at (1.2, 0.6);    
\coordinate (B) at (3.8, 3.8);    
\coordinate (C) at (3.4, 1.1);   
\coordinate (D) at (5.4, 0.1);    
\coordinate (E) at (7.9, 3.3);    
\coordinate (F) at (6.2, 2.05);    
\coordinate (G) at (6.6, 5.0);    
\coordinate (H) at (7.5, -0.5);   

\draw[black, thick]
  (A) -- (B) -- (C) -- cycle
 
  (B) -- (D)
  (B) -- (E)
  (B) -- (F)
  (B) -- (G)
  (C) -- (D)
  (D) -- (F)
  (E) -- (F)
  (E) -- (G)
  (E) -- (H)
  (F) -- (H);

\draw[teal!80!green, thick]
 (A) -- (6.5,4.5);
 
\draw[red!80!black, thick]
 (A) -- (7.8,2.9);

\draw[teal!80!blue, thick]
 (A) -- (E);
 
\end{tikzpicture}
    \caption{Unfoldings after four reflections of two billiard orbits (in green and red), starting at the same vertex. The first three steps of unfoldings are the same for both trajectories, but the fourth differs because the shortest generalized diagonal belonging to the angular segment between the trajectories (in blue) has combinatorial length three. This generalized diagonal defines a cutting point in the partition $\xi_n$ of the initial angle for all $n\geq 3.$}
    \label{unfolding}
\end{figure}
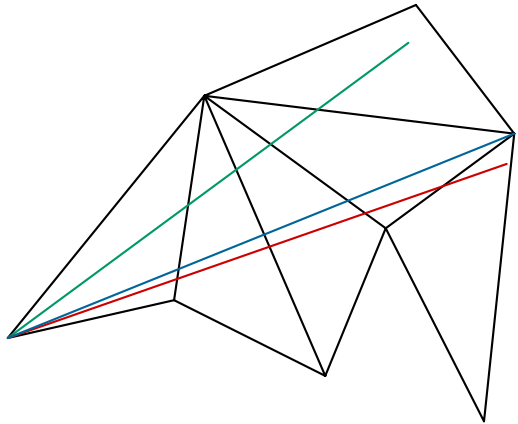

The classical construction to study polygonal billiards is an unfolding process introduced independently by R. Fox and R. Kershner \cite{FK} and 40 years later by A. Zemljakov and A. Katok \cite{KZ}. Instead of reflecting a billiard trajectory on a side of the polygon $Q$, we reflect $Q$ on the side and unfold the trajectory into a straight line. This process is repeated at each reflection: each successive copy of the polygon is obtained from the previous one by reflecting it across the side encountered by the straightened trajectory. See Figure \ref{unfolding} for an illustration.

Fix a vertex of the triangle and consider the corresponding angular segment $\theta$. Generalized diagonals of combinatorial length at most $n$, starting at this vertex, divide $\theta$ into smaller angular segments. Identifying the initial angular segment $\theta$ with an interval, we obtain an interval partition, denoted $\xi_n$.

The criterion for weakly exponential growth for a given triangular billiard is that two generalized diagonals starting or terminating at the same vertex cannot be arbitrarily close. It can be formulated in terms of interval partitions introduced above. The criterion comes from \cite{Sch2}, although it is not explicitly stated. We repeat the suite of lemmas from \cite{Sch2} that leads to the criterion. The weakly exponential growth proof (presented in subsections 4.1 and 4.2 of \cite{Sch2}) works for all triangles satisfying a property from Lemma 4.1 (used for a construction in Lemma 4.7 of \cite{Sch2}). We single out this property as a criterion for weakly exponential growth. 

\begin{thm}[The criterion for weakly exponential growth, \cite{Sch2}]\label{criterion}
A triangular billiard is weakly exponential if there is a universal constant $a>0$ with the following property:
    
    Any interval $I$ of the $n$-th partition $\xi_n$ (corresponding to any vertex of the
triangle) is of the length $|I|>e^{-an^2}$.
\end{thm}

We show that this criterion is satisfied for all triangles whose two angles are trigonometrically algebraic. In the general case, the criterion is satisfied for almost all triangles in the Lebesgue-measure sense (parameterized by two smaller angles $\alpha$, $\beta$; see \cite{Sch2} for the proof). This yields a proof of weakly exponential growth for almost every triangle.

\begin{thm}[Weakly exponential growth, \cite{Sch2}]
    The growth of the combinatorial complexity function of almost all triangular billiards is weakly exponential.
\end{thm}

The length of the interval $I$ of the $n$-th partition $\xi_n$ is determined by the coordinates of the endpoints of the unfoldings of the generalized diagonals bounding $I$. The coordinates of an unfolding of a generalized diagonal of length at most $k\le n$ are defined by a sequence of $k$ reflections across the sides or, equivalently, a sequence of $k/2$ rotations by angles $2\alpha$, $2\beta$, and $\pi-\alpha-\beta$. So, the length of $I$ can be estimated using an integer polynomial in the variables $\sin\alpha,\cos\alpha,\sin\beta,\cos\beta$ with bounded total degree and coefficients. 

\begin{lem}[Formulated in the proof of Lemma 4.1 of \cite{Sch2}]\label{length_estim}
Let $I\in\xi_n$ be an interval of the $n$-th partition corresponding to any vertex of the
triangle. Then
$$|I| > \frac{\ell\lvert M_{I,n}(\alpha, \beta)\rvert}{n^2},$$
where $\ell>0$ is a universal constant depending only on $\alpha,\beta$ and
\[M_{I,n}(\alpha, \beta) := M_{I,n}(\sin \alpha, \cos \alpha, \sin \beta, \cos \beta) =\]
\[=\sum_{\substack{k_1,k_2,k_3,k_4 \\ k_1+\dots+k_4 \le 4n}}
a_{k_1k_2k_3k_4}
\sin^{k_1}\alpha \cos^{k_2}\alpha
\sin^{k_3}\beta \cos^{k_4}\beta\]
is a nonzero integer polynomial depending on the number of partition $n$ and a chosen interval $I$. Moreover, the total degree and the coefficients of $M_{I,n}$ are bounded linearly and exponentially, respectively, on the number $n$ of the partition $\xi_n$. More precisely, $\deg M_{I,n} \le 4n$ (i.e. $k_1+k_2+k_3+k_4\le4n$) and $ \max|a_{k_1, k_2, k_3, k_4}| \le n^{14} 4^{2n+6}$.

\begin{rmrk}
    The precise form of the bound on the coefficients (respectively, on the degree) is not important; we only use the fact that it depends exponentially (respectively, linearly) on $n$, where $n$ is the index of the partition to which $I$ belongs. 

\end{rmrk}

\begin{rmrk}
    Geometrically, the polynomial $M_{I,n}$ represents the area of a parallelogram spanned by unfoldings of generalized diagonals of combinatorial length at most $n$ (namely, those bounding $I$; see Lemma 4.1 in \cite{Sch2} for details). We denote by $\mathcal{F}_n$ the set of all such polynomials; all of them satisfy the bounds on degree and coefficients from the previous lemma. Since these polynomials are defined by two sequences of at most $n/2$ rotations on angles $2\alpha$, $2\beta$ and $2(\pi-\alpha-\beta)$, the cardinality of $\mathcal{F}_n$ satisfies $|\mathcal{F}_n|<e^{tn}$ for some $t>0$ (depending only on $\alpha,\beta$). 
\end{rmrk}
\end{lem}

An elementary corollary of the criterion for weakly exponential growth is the following statement, providing the first explicit example of triangular billiards of weakly exponential growth. 

\begin{cor}\label{rational}
     Assume that the angles $\alpha,\beta$ are trigonometrically rational. Then the growth of the combinatorial complexity function of a billiard in the triangle $(\alpha,\beta,\pi-\alpha-\beta)$ is weakly exponential. 
\end{cor}

\begin{rmrk}
    Clearly, the family of such triangles is infinite. A triangle of this family can be viewed as a couple of Pythagorean triples, defining $\sin\alpha,\cos\alpha,\sin\beta,\cos\beta$. 
\end{rmrk}

\begin{Proof}
    Since $|M_{I,n}(\alpha, \beta)|\neq 0$ and $M_{I,n}$ is a polynomial with integer coefficients of degree at most $4n$, we obtain a lower bound for $|M_{I,n}(\alpha, \beta)|$. Let  $\sin\alpha=\frac{p_1}{q_1}$, $\cos\alpha=\frac{p_2}{q_2}$, $\sin\beta=\frac{p_3}{q_3}$, $\cos\beta=\frac{p_4}{q_4}$. Substituting variables, we have
    \[M_{I,n}(\sin\alpha, \cos\alpha, \sin\beta, \cos\beta)=\frac{M_1(p_i,q_i)}{q_1^{4n} q_2^{4n} q_3^{4n} q_4^{4n}},\] where $M_1(p_i,q_i)$ is again a nonzero polynomial with integer coefficients. Thus, $|M_1(p_i,q_i)|\geq 1$ and \[|M_{I,n}(\sin\alpha, \cos\alpha, \sin\beta, \cos\beta)|\geq \frac{1}{(q_1q_2q_3q_4)^{4n}}.\] 
    By Lemma \ref{length_estim} for any interval $I$ of the $n$-th partition corresponding to any vertex of the
triangle
    \[|I|>\frac{\ell|M_{I,n}(\alpha,\beta)|}{n^2}\geq \frac{\ell}{n^2(q_1q_2q_3q_4)^{4n}}\] and there is a universal constant $a$ such that $|I|>e^{-an}>e^{-an^2}$ for any $n$. The criterion for weakly exponential growth is satisfied; therefore, the growth of a billiard in such a triangle is weakly exponential. 
\end{Proof}

Here we used the easiest version of Liouville’s inequality, that if the value of an integer polynomial at an integer number is nonzero, then its absolute value is at least 1. We proove the same result for trigonometrically algebraic angles, using the general form of the Liouville inequality. In general, less ``convenient'' the values of variable in an integer polynomial are, more conditions on a polynomial and more general form of Liouville inequality we should use.
\begin{defn}
   For a polynomial $P(x_1,\dots,x_t)$ the \textit{height} $H(P)$ of $P(x_1,\dots,x_t)$ is the maximum of  the absolute values of its coefficients, and the \textit{length} $L(P)$ the sum of the absolute values of its coefficients.  
\end{defn}
Note that for the polynomial in $t$ variables $L(P)\leq (\deg P+1)^t\cdot H(P)$, where the notation $\deg P$ is used for the total degree of $P$, i.e. $\deg(\sum_{q_1\dots q_t} a_{q_1\dots q_t} x_1^{q_1}\cdots x_t^{q_t})=\max_{q_1\dots q_t}(q_1+\dots+q_t)$.

We need the special case of the following theorem.

\begin{thm}[Liouville inequality, Theorem B.10 in \cite{B}.]
Let $f(x_1,\dots,x_t)$ be a nonzero polynomial in $t$ variables with integer coefficients.
Let $\gamma_1,\dots,\gamma_t$ be algebraic numbers in a number field $K$ of degree $d$, such that $f(\gamma_1,\dots,\gamma_t)\neq 0$.

Then for every place $v$ of $K$ the following estimate holds:
\[
\log |f(\gamma_1,\dots,\gamma_t)|_v
\ge
-\frac{d}{d_v}
\left(
\log L(f) + \sum_{i=1}^t \deg_{x_i}(f)\, h(\gamma_i)
\right),
\]
where \( d_v \) denotes the local degree at the place \( v \) and $h(\gamma_i)$ denotes the (logarithmic) Weil height of $\gamma_i$.

\end{thm}

We will only use this theorem in the case of real algebraic $\gamma_1,\dots,\gamma_t$. In this case each $h(\gamma_i)$ is a non-negative real constant, depending only on $\gamma_i$, and $d$ is a positive integer, depending only on $(\gamma_1,\dots,\gamma_t)$. For details and definitions in general case see \cite{B}. In the case of real algebraic $\gamma_1,\dots,\gamma_t$ the result could be rewritten in the following form.

\begin{lem}\label{Liouville} Let $f(x_1,\dots,x_t)$ be a polynomial in $t$ variables with integer coefficients and $\gamma_1,\dots,\gamma_t$ be real algebraic numbers such that $f(\gamma_1,\dots,\gamma_t)\neq 0$.

Then the following estimate holds:
    \[
\lvert f(\gamma_1\dots\gamma_t)\rvert
\ge
(\deg f+1)^{-td}\,(H(f))^{-d}\,(c_{\gamma_1\dots\gamma_t})^{-\deg f},
\]
where constants $d\in \mathbb{N}$ and \( c_{\gamma_1\dots\gamma_t}>0 \) depend only on \( \gamma_1,\dots,\gamma_t \).
\end{lem}
\begin{Proof}
%This lemma is the special case of the previous theorem. 
We simply rewrite the Liouville inequality in the case of $t$ real algebraic $\gamma_i$. Let $K=\mathbb{Q}(\gamma_1,\dots,\gamma_t)$. For an Archimedean place \( v \) corresponding to the standard embedding $K \hookrightarrow \mathbb{C}, \quad x \mapsto x$ one has $|x|_v = |x|$ and $ d_v = 1$. Denote by $\deg f$ the degree of $f$ with respect to all variables; $\deg_{x_i} f \le \deg f$. Using $L(f) \le H(f)\,(\deg f+1)^t$, we obtain the estimate
\[
\lvert f(\gamma_1\dots\gamma_t)\rvert
\ge
(L(f))^{-d}\, e^{-d\sum_{i=1}^t \deg_{x_i}(f)\, h(\gamma_i)}\;\ge\]
\[
\ge\;
((\deg f+1)^t\,H(f))^{-d}\, e^{-d\sum_{i=1}^t \deg(f) h(\gamma_i)}
\;\ge\;\]
\[\ge \;
(\deg f+1)^{-td}\,(H(f))^{-d}\,(c_{\gamma_1\dots\gamma_t})^{-\deg f}\; ,
\]
where $d$ and \( c_{\gamma_1\dots\gamma_t} \) depend only on \( \gamma_1,\dots,\gamma_t \).
\end{Proof}

This yields the following theorem, which explicitly describes a vast family of triangular billiards with weakly exponential growth.

\begin{thm}
 Assume that the angles $\alpha,\beta$ are trigonometrically algebraic. Then the growth of the combinatorial complexity function of a billiard in the triangle $(\alpha,\beta,\pi-\alpha-\beta)$ is weakly exponential.   
\end{thm}
\begin{Proof}

For an interval \( I \in \xi_n \),
\[
|I|
\ge
\frac{\ell}{n^2}\,
\bigl| M_{I,n}(\sin\alpha,\cos\alpha,\sin\beta,\cos\beta) \bigr|,
\]
where the conditions on $M_{I,n}$ are described in \ref{length_estim}. 

By Liouville's inequality in the form of Lemma \ref{Liouville},
\[
|M_{I,n}(\sin\alpha,\cos\alpha,\sin\beta,\cos\beta)|
\ge
(\deg M_{I,n} + 1)^{-4d}\,
(H(M_{I,n}))^{-d}\,
(c_{\alpha,\beta})^{-\deg M_{I,n}},
\]
where \( c_{\alpha,\beta} \) is a constant depending only on \( \alpha,\beta \) and \( d \) is the degree of the number field $K = \mathbb{Q}(\sin\alpha,\cos\alpha,\sin\beta,\cos\beta)$ depending, again, only on \( \alpha\) and \(\beta \). 
Therefore,
\[
|M_{I,n}(\sin\alpha,\cos\alpha,\sin\beta,\cos\beta)|
\ge
(4n+1)^{-4d}\,
\bigl( n^{14}\,4^{\,2n+2} \bigr)^{-d}\,
(c_{\alpha,\beta})^{-4n}.
\]
The estimate implies that there is a constant \( a_0 > 0 \) such that \(|M_{I,n}| \ge e^{-a_0 n},\) therefore there is a constant \( a > 0 \) such that $|I|>e^{-an}>e^{-an^2}$ and the criterion for weakly exponential growth is satisfied.
\end{Proof}

\begin{rmrk}
 Note that in the case of trigonometrically algebraic angles we obtain an estimate that is much stronger than required by the criterion: we have $e^{-an}$ instead of $e^{-an^2}$. Nevertheless, it does not allow us to improve the complexity estimate, following the strategy of D. Scheglov. Indeed, for $k\in[1,2]$, replacing $e^{-an^2}$ with $e^{-an^k}$ in the criterion changes only the negligible polynomial part of the complexity estimate (for more details, see Lemma 4.7 in \cite{Sch2}). In the case of trigonometrically algebraic angles, one obtains a better estimate, similarly to the proof of the strong Diophantine condition for elements of $SU(2)$ with algebraic entries (A. Gamburd, D. Jakobson and P. Sarnak, \cite{GJS}). Since the question of the strong Diophantine condition is still open, to obtain ``almost everywhere'' result the criterion is stated using  $e^{-an^2}$. It corresponds to the weak Diophantine condition established by V. Kaloshin and I. Rodnianski in \cite{KR}. Note also that in the following section, dedicated to ``almost all'' results, we obtain the weak Diophantine condition: $|I|>e^{-an^2}$.
\end{rmrk}

\section{Distribution of angles}\label{sec3}

Now we are interested in the distribution of angles that satisfy the criterion for weakly exponential growth. Is it true that for almost every triangle with one fixed angle the combinatorial complexity is weakly exponential? We first consider a simple particular case, illustrating the strategy presented in \cite{Sch2}. We need the following result.

\begin{thm}[Proposition 1 in \cite{KR}\label{KalRod} reduced to the form of Corollary 4.1 in \cite{Sch2}] 
Set three conditions on a trigonometric polynomial $P(\alpha):=P(\sin\alpha,\cos\alpha)$:
\begin{enumerate}
         \item $P(\sin\alpha,\cos\alpha)$ is a nonzero polynomial with integer coefficients;
         \item $\deg P\leq m$;
         \item $H(P)\leq(m2^m)^2$.
     \end{enumerate}
     
There exist universal constants $R, h>0$ such that for any polynomial $P(\alpha)$ satisfying conditions 1-3 the following estimate holds:
\[\Leb\SizedSet{\alpha\given |P(\sin\alpha, \cos\alpha)|<e^{-Rm^{2}}} <e^{-Rm^{2}}\rbrace<e^{-hm}.\]
\end{thm}
Here and throughout, $\Leb$ denotes the Lebesgue measure.
\begin{rmrk}
    Recall that $\mathcal{F}_n$ denotes the set of polynomials $M_{I,n}$ representing the area of the parallelogram spanned by two neighboring generalized diagonals of length at most $n$; the cardinality of $\mathcal{F}_n$ satisfies $|\mathcal{F}_n|<e^{tn}$ for some $t>0$. 
\end{rmrk}

The first simple example is a triangle with a fixed angle $\pi/2$.
\begin{lem}\label{right}
    The complexity function of almost every right triangular billiard is weakly exponential.
\end{lem}

\begin{Proof}
For a right triangular table $\beta = \frac{\pi}{2} - \alpha$, $\sin \beta = \cos \alpha$, $\cos \beta = \sin \alpha$. Then the length of interval $I$ of a partition $\xi_n$ satisfies
\[
|I|>\frac{\ell|M_{I,n}(\alpha,\beta)|}{n^2}
= \frac{\ell|M_{I,n}(\alpha, \frac{\pi}{2}-\alpha)|}{n^2},\quad\text{with}
\]
\[
M_{I,n}(\alpha, \frac{\pi}{2}-\alpha) = \sum_{\substack{k_1 , k_2, \\ k_3, k_4}} a_{k_1k_2k_3k_4}
\sin^{k_1} \alpha \cos^{k_2} \alpha
\sin^{k_3}\left(\frac{\pi}{2}-\alpha\right)
\cos^{k_4}\left(\frac{\pi}{2}-\alpha\right)=
\]
\[
= \sum_{\substack{k_1 , k_2, \\ k_3, k_4}} a_{k_1k_2k_3k_4}
\sin^{k_1 + k_4} \alpha \cos^{k_2 + k_3} \alpha,
\]
where $\ell>0$ is a universal constant, $a_{k_1k_2k_3k_4} \le n^{14} 4^{2n+6}$, $a_{k_1k_2k_3k_4} \in \mathbb{Z}$, and $k_1 + k_2 + k_3 + k_4 \le 4n$. Rewriting,
\[
M_{I,n}(\alpha,\frac{\pi}{2}-\alpha)=\sum_{i+j\le 4n} b_{ij}\sin^i\alpha\cos^j\alpha, \text{ where  } b_{ij}=\sum_{\substack{k_1+k_4=i, \\k_2+k_3=j}} a_{k_1k_2k_3k_4}.\] It follows that $|b_{ij}| \le (i+1)(j+1)\max |a_{k_1k_2k_3k_4}|
\le (4n+1)^2 n^{14} 4^{2n+6} \le 2^{kn}$ for some constant $k$, independent of $n$. By Lemma \ref{length_estim}, all polynomials of $\mathcal{F}_n$ satisfy the conditions of Theorem \ref{KalRod} for $m$ sufficiently large (in particular, for all $m\ge\max\{4n,kn\}$).

To complete the proof, we repeat the argument of \cite{Sch2}. Let $m = F n$, where $F\in\mathbb{N}$ is a constant to be chosen later. By Theorem \ref{KalRod} there are universal constants $ R, h > 0$ such that for any polynomial $M_{I,n}$ of degree at most $m$ the estimate
\[
\Leb\SizedSet{\alpha\given |M_{I,n}(\alpha)| < e^{-R m^2}} < e^{-hm}
\] holds. In particular, when $m\ge 4n$ it holds for any polynomial $M_{I,n}\in\mathcal{F}_n$.

Define the set of angles
${\mathcal{B}}_n=\SizedSet{\alpha\given\exists M_{I,n}\in{\mathcal{F}}_n \text{ such that }|M_{I,n}(\alpha)|<e^{-Rm^2}}$. As $|\mathcal{F}_n|<e^{tn}$ for some $t>0$ then $\Leb({\mathcal{B}}_n)<|\mathcal{F}_n|\cdot e^{-hm}=e^{tn-hm}=e^{(t-hF)n}$. We choose the constant $F>\max\{4,k,t/h\}$ to obtain $\sum \Leb({\mathcal{B}}_n)<\infty$. By the Borel-Cantelli lemma for Lebesgue-almost any
$\alpha$ and all sufficiently large $n$ for any polynomial $M_{I,n}\in\mathcal{F}_n$ we have $|M_{I,n}(\alpha)|\geq e^{-RF^2n^2}$. Choosing $a> RF^2$ (depending on $\alpha$) we obtain that for almost any $\alpha$ the length of interval $|I|>e^{-an^2}$ and the growth of the complexity function is weakly exponential.

\end{Proof}

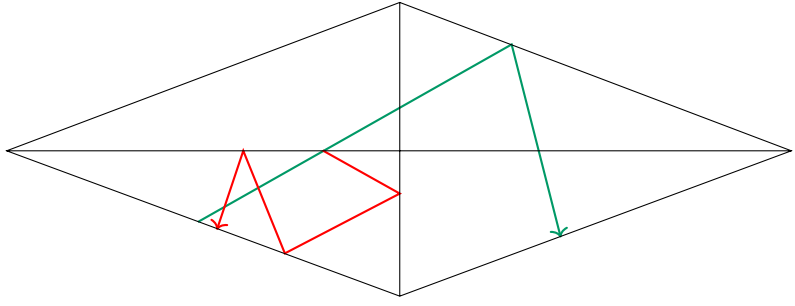
\begin{figure}
\begin{tikzpicture}
\draw[teal!80!green, thick, ->]
    (2.55,-1.945) -- (6.7,0.4) -- (7.35,-2.15);

\draw[red, thick, ->]
    (4.2,-1) -- (5.22,-1.57) -- (3.7,-2.36) -- (3.15,-1.005) -- (2.8,-2.05);

\draw 
(0.02,-1.005)  -- (5.22,-2.925) -- (10.4,-1.005) -- (5.22,0.955) -- cycle
(0.02,-1.005)  -- (10.4,-1.005)
(5.22,-2.925) -- (5.22,0.955);

\end{tikzpicture}
    \caption{A trajectory of billiard ball in a rhombus (green) and the corresponding trajectory in the right triangle (red). The trajectory in rhombus has two links, in the triangle -- five links; in general, if the trajectory in rhombus has $n$ links, then it has at most $3n$ links in the triangle, providing $N_{\text{rhombus}}(n)\leq N_{\text{right triangle}}(3n)$.}
    \label{rhombus}
\end{figure}

Note that from weakly exponential growth for almost all right triangles it follows weakly exponential growth for almost all isosceles triangles and almost all rhombus. Indeed, by elementary geometric reasoning $N_{\text{rhombus}}(n)\leq N_{\text{right triangle}}(3n)$, where a right triangle corresponding to rhombus is one of four equal right triangles obtained from rhombus by dividing it by diagonals. Every generalized diagonal of combinatorial length $n$ in a rhombus is a generalized diagonal in the corresponding right triangle of combinatorial length at most $3n$ (see Figure \ref{rhombus} and \cite{Sch1}, Lemma 2.2 for more details). By the same reasoning, $N_{\text{isosceles triangle}}(n)\leq N_{\text{right triangle}}(2n)$ where a right triangle corresponding to an isosceles triangle is one of two right triangles obtained from the original  by dividing it by the altitude to the base. Similarly, every rhombus with trigonometrically algebraic angles has weakly exponential complexity growth.

\begin{rmrk}
    Using the same strategy, one can also treat the case of the $(k_1\alpha,k_2\alpha,\pi-(k_1+k_2)\alpha)$ triangle for any fixed $k_1,k_2\in\mathbb{N}$. Indeed, using trigonometric identities, we can represent nonzero $M_{I,n}(k_1\alpha,k_2\alpha)$ as $M_1(\alpha)$ with integer coefficients bounded exponentially in $n$ and degree bounded linearly in $n$. Applying the Borel-Cantelli lemma, we obtain that for any $k_1,k_2\in\mathbb{N}$ almost all $(k_1\alpha,k_2\alpha,\pi-(k_1+k_2)\alpha)$-triangular billiards has a complexity function of weakly exponential growth. Equivalently, for any  $k\in\mathbb{Q}$ and almost any $\alpha$ the billiard in the triangle $(\alpha,k\alpha,\pi-(k+1)\alpha)$ is of weakly exponential complexity.
\end{rmrk}

For a fixed angle $\alpha$ with algebraic but non-rational values of trigonometric functions, the trigonometric polynomial $M_{I,n}(\beta):=M_{I,n}(\alpha,\beta)$ can no longer be normalized to the form of an integer polynomial, and thus Theorem \ref{KalRod} is not applicable. We prove another technical lemma of the same form to treat the case of a fixed algebraic angle. We need the following important lemma.

\begin{lem}[Lemma 4.1 in \cite{DM} reduced to the form of Lemma 5 in \cite{KR}]
    Let $P(x)$ be a polynomial of degree $\leq n$. Denote by $\norm{P}=\max_{x\in[-1,1]}|P(x)|$. Then
\[\Leb\SizedSet{x\in[-1,1] \given |P(x)|< e^{-Rn^2}}< 4n(n+1)^{1/n}\left(\frac{e^{-Rn^2}}{\lVert P\rVert}\right)^{1/n}\]
\end{lem}

For monic polynomials, this result takes the following form.

\begin{lem}\label{unitary}
    Let $P(x)$ be a monic polynomial of degree $\leq n$. Then for any $R>0$ there is $h>0$ such that
\[\Leb\SizedSet{x\in[-1,1]\given \lvert P(x)\rvert < e^{-Rn^2}} < e^{-hn}\]
\end{lem}

\begin{Proof}
    Let $\deg P = k\leq n$. The $k$-th normalized Chebyshev polynomial of the first kind minimizes the supremum norm among monic polynomials of degree $k$ in $[-1,1]$. So for $T_k(x)=2^{k-1}x^k+\dots$ one has
\[
    \frac{1}{2^{k-1}}\lVert T_k(x)\rVert=\left\lVert\frac{1}{2^{k-1}}T_k(x)\right\rVert\leq\lVert P(x)\rVert;
\]
\[
\frac{1}{\lVert P(x)\rVert}\leq \frac{2^{k-1}}{\lVert T_k(x)\rVert}.
\]
Moreover, for every $k$ the norm of the non-normalized Chebyshev polynomial $\lVert T_k(x)\rVert=1$, which implies
\[
\frac{1}{\lVert P(x)\rVert}\leq 2^{k-1}\le 2^{n-1}.
\]
For any chosen $R$ the estimate from the previous lemma takes form
\[ \Leb\SizedSet{x\in[-1,1]\given \lvert P(x)\rvert < e^{-Rn^2}} < 4n(n+1)^{1/n}\left(\frac{e^{-Rn^2}}{\norm{P}}\right)^{1/n}\leq
\]
\[
\leq 4n(n+1)^{1/n}\left(2^{n-1}e^{-Rn^2}\right)^{1/n}= 4n(n+1)^{1/n}2^{\frac{n-1}{n}}e^{-Rn},
\]

therefore, there is a constant $h>0$ such that
\[\Leb\SizedSet{x\in[-1,1]\given \lvert P(x)\rvert < e^{-Rn^2}} < e^{-hn}. \]
\end{Proof}

\begin{thm}
    Let $\alpha$ be a fixed trigonometrically algebraic angle. Then for almost all $\beta$, the complexity function of the billiard in the triangle $(\alpha,\beta,\pi-\alpha-\beta)$ is weakly exponential.
\end{thm} 

\begin{Proof}
We show that in the case of fixed trigonometrically algebraic $\alpha$, the criterion  is satisfied for almost every $\beta$.

\[
M_{I,n}(\sin\alpha, \cos\alpha, \sin\beta, \cos\beta)
=
\sum_{\substack{k_1 , k_2, \\ k_3, k_4}}
a_{k_1k_2k_3k_4}
\sin^{k_1}\alpha \cos^{k_2}\alpha
\sin^{k_3}\beta \cos^{k_4}\beta
\]
To reduce the trigonometric polynomial in two dependent variables $\sin\beta$, $\cos\beta$ to a polynomial in one variable, let us apply the following change of variables:
\[\sin \beta = \frac{2t}{t^2 + 1},\quad
\cos \beta = \frac{1 - t^2}{t^2 + 1} , \quad\text{where } t = \tan \frac{\beta}{2}. \]

The consider the case of $\beta \le \frac{\pi}{2}$ and do not treat the case $\beta \ge \frac{\pi}{2}$ in detail -- the proof is the same, but the change of variables is $t := \cotan \frac{\beta}{2}$  and not $\tan \frac{\beta}{2}$. The reason for distinguishing these cases is that for $\beta \le \frac{\pi}{2}$ one has $t=\tan \frac{\beta}{2} \in (0, 1]$ (respectively, for $\beta \in (\frac{\pi}{2},\pi]$ one has $t=\cotan  \frac{\beta}{2} \in (0, 1]$) which is used below. 

After the change of variables
\[
M_\alpha(t):=M_{I,n}\left(\sin\alpha,\cos\alpha, \frac{2t}{t^2 + 1}, \frac{1 - t^2}{t^2 + 1}\right) =\]
\[= \sum_{\substack{k_1 , k_2, \\ k_3, k_4}} a_{k_1k_2k_3k_4}
\sin^{k_1}\alpha \cos^{k_2}\alpha \left(\frac{2t}{t^2+1}\right)^{k_3}
\left(\frac{1-t^2}{t^2+1}\right)^{k_4};
\]
\[
(t^2+1)^{4n} M_\alpha(t)
= \sum_{\substack{k_1 , k_2, \\ k_3, k_4}} a_{k_1k_2k_3k_4}
\sin^{k_1}\alpha \cos^{k_2}\alpha (2t)^{k_3} (1-t^2)^{k_4} (t^2+1)^{4n-k_3-k_4}=
\]
\[
= \sum_{\substack{ k_1 , k_2 \\ k_3, k_4}} \sum_{\substack{m\le k_4\\ l\le 4n-k_3-k_4}}
2^{k_3} a_{k_1\dots k_4}
\sin^{k_1}\alpha \cos^{k_2}\alpha
(-1)^m \binom{k_4}{m}
\binom{4n-k_3-k_4}{l}
t^{k_3 + 2m + 2l}=
\]
\[
= \sum_{i \le 8n} b_i(\alpha) t^i,
\]
where
\[
b_i(\alpha) =
\sum_{\substack{k_1 , k_2, \\ k_3, k_4}}\quad
\sum_{k_3 + 2m + 2l = i}
2^{k_3} a_{k_1k_2k_3k_4}
\sin^{k_1}\alpha \cos^{k_2}\alpha
(-1)^m \binom{k_4}{m}
\binom{4n-k_3-k_4}{l}
\]
and at least one of $b_i(\alpha)$ is nonzero. Note that the degree of $(t^2+1)^{4n} M_\alpha(t)$ is linearly bounded on $n$: $\deg\left((t^2+1)^{4n} M_\alpha(t)\right)\le12n$. To normalize $M_\alpha(t)$ to a monic polynomial by ``bounded multiplication'', we would like to obtain a lower bound for the first nonzero coefficient $b_i(\alpha)$. Write
\[
b_i(\alpha) = \sum_{k_1k_2} c_{k_1k_2}
\sin^{k_1}\alpha \cos^{k_2}\alpha,
\]
where
\[
c_{k_1k_2} =
\sum_{k_3,k_4}
2^{k_3} a_{k_1k_2k_3k_4}
\sum_{k_3+2m+2l=i}
(-1)^m \binom{k_4}{m}
\binom{4n-k_3-k_4}{l};
\]
\[
|c_{k_1k_2}|
\le \sum_{k_3,k_4}
2^{k_3} \max |a_{k_1k_2k_3k_4}|
\sum_{m,l}
\binom{k_4}{m}
\binom{4n-k_3-k_4}{l}\le
\]
\[
\le 2^{4n} (4n+1)^2 \max |a_{k_1k_2k_3k_4}|
\le e^{a_0 n}
\]
for some \(a_0 > 0\).

By the Liouville inequality of the form \ref{Liouville} for all nonzero $b_i(\alpha)$ and $d=[\mathbb{Q}(\sin\alpha,\cos\alpha):\mathbb{Q}]$ the following estimate holds:
\[
|b_i(\alpha)| \ge (4n+1)^{-4d}(e^{a_0n})^{-d}(c_\alpha)^{-4n}\ge e^{-a_1 n}
\]
for some \(a_1 > 0\). Therefore
\[
(t^2+1)^{4n} M_\alpha(t) = \sum_i b_i(\alpha) t^i,
\]
where at least one of $b_i(\alpha)$ is nonzero and all nonzero $b_i(\alpha)$ satisfy $|b_i(\alpha)| \ge e^{-a_1 n}$. Choose the largest index $i$ such that $b_i(\alpha)\neq 0$. Dividing by $b_i(\alpha)$ gives a monic polynomial $\frac{(t^2+1)^{4n} M_{I,n}(t)}{b_i(\alpha)}$ of degree at most $16n$.

We show that for every $R>0$ there is a constant $h_1>0$ such that for any polynomial $M_{I,n}\in\mathcal{F}_n$ one has the estimate $\Leb\SizedSet{t \given \abs{M_{I,n}(t)} < e^{-R n^2}}< e^{-h_1 n}.$ If $\left\lvert M_\alpha(t)\right\rvert< e^{-Rn^2}$ then, as the value of $t$ is bounded, we obtain $\left\lvert(t^2+1)^{4n}M_\alpha(t)\right\rvert<e^{-Rn^2+a_0n}$ for some constant $a_0>0$ and
\[\left\lvert\frac{(t^2+1)^{4n}M_\alpha(t)}{b_i(\alpha)}\right\rvert< e^{-Rn^2+a_0n+a_1n},
\]
so for some $R_1>0$ and all sufficiently large $n$ we have the estimate
\[\left\lvert\frac{(t^2+1)^{4n}M_\alpha(t)}{b_i(\alpha)}\right\rvert< e^{-R_1n^2}.
\]
Therefore, for all $R>0$ there is a constant $R_1>0$ such that: \[
\SizedSet{t \given \lvert M_{I,n}(t)\rvert < e^{-R n^2}}
\;\subseteq\;
\SizedSet{t \given
\left\lvert \frac{(t^2+1)^{4n} M_{I,n}(t)}{b_i(\alpha)}\right\rvert  < e^{-R_1 n^2}}.
\]
The polynomial $\frac{(t^2+1)^{4n} M_{I,n}(t)}{b_i(\alpha)}$ is monic of degree at most $16n$, so by Lemma \ref{unitary} there is a constant $h>0$ such that
\[
\Leb\SizedSet{t \given
\left\lvert\frac{(t^2+1)^{4n} M_{I,n}(t)}{b_i(\alpha)}\right\rvert < e^{-R_1n^2} = e^{-\frac{R_1}{16^2} (16n)^2}}
\;<\;
e^{-h n}.
\]
So $\Leb\SizedSet{t \given M_{I,n}(t) < e^{-R n^2}}< e^{-h n}$ and hence there is a constant $h_1>0$ such that
\[\Leb\SizedSet{\beta \given |M_{I,n}(\sin\beta,\cos\beta)| < e^{-R n^2}} <
e^{-h_1 n}.\]
Repeating the argument for $\beta\in(\pi/2,\pi]$, we obtain that for any $R>0$ there is $h_1 > 0$ such that
\[
\Leb\SizedSet{\beta\in(0,\pi]\given |M_{I,n}(\beta)| < e^{-R n^2}} < e^{-h_1n}.
\]

This estimate holds for all $M_{I,n}\in\mathcal{F}_n$ and the constant $h_1$ is independent of $M_{I,n}$. Note that for any $m\ge n$ we have $\mathcal{F}_n\subseteq \mathcal{F}_m$, so we can reformulate the result in a more general form: for any $R>0$ there is $h_1 > 0$ such that for any polynomial $M_{I,n}\in\mathcal{F}_n$ satisfying the conditions of Lemma \ref{length_estim} and all $m\geq n$, the estimate
\[
\Leb\SizedSet{\beta\in(0,\pi]\given |M_{I,n}(\beta)| < e^{-R m^2}} < e^{-h_1m}
\] holds. 

To complete the proof, we use the same Borel-Cantelli argument as in Lemma \ref{right}. Let $m = F n$, where $F\in\mathbb{N}$ is a constant to be chosen later. We fix $R$ and define the set of angles
\[{\mathcal{B}}_n=\SizedSet{\beta\given\text{ there is a polynomial } M_{I,n}\in{\mathcal{F}}_n \text{ such that }|M_{I,n}(\beta)|<e^{-Rm^2}}.\]

As $|\mathcal{F}_n|< e^{tn}$ for some $t>0$ then $\Leb({\mathcal{B}}_n)<|\mathcal{F}_n|\cdot e^{-h_1m}=e^{tn-h_1m}=e^{(t-h_1F)n}$. We choose the constant $F>t/h_1$ to obtain $\sum \Leb({\mathcal{B}}_n)<\infty$. By the Borel-Cantelli lemma for Lebesgue-almost any
$\beta$ and all sufficiently large $n$ for any polynomial $M_{I,n}\in\mathcal{F}_n$ we have $|M_{I,n}(\alpha,\beta)|\geq e^{-RF^2n^2}$. Hence, when $\alpha$ is a fixed trigonometrically algebraic angle, for almost any $\beta$ there exists a constant $a>0$ (depending on $\alpha$) such that the length of the interval $|I|>e^{-an^2}$. Therefore, the criterion is satisfied, and the growth of the complexity function is weakly exponential.
\end{Proof}

\printbibliography
\end{document}